\documentclass[12pt]{amsart}
\usepackage[all]{xy}
\usepackage{graphicx}
\usepackage{amssymb}
\makeatletter
\newdimen\sectionruledimen
\def\makeparrule{%
 \def\par{%
  \endgraf\nobreak\vskip\lineskip\nointerlineskip
  \hbox to\hsize{\hskip\sectionbackskip\leaders\hrule height
  \sectionruledimen\hfil}%
  }%
 }%

\def\section{%
 \@startsection {section}{1}{\sectionbackskip}{-10pt plus -1.2ex minus -.2ex}%
  {.5pt}{\normalsize\bf\makeparrule}%
 }%
\newdimen\sectionbackskip
\sectionbackskip\z@
\newdimen\sectionruledimen
\sectionruledimen\z@
\makeatother

\newenvironment{theor}{\smallskip\begin{trivlist}
 \item[\hspace{\labelsep}{\noindent\bf Theorem.}]\it
 }{\end{trivlist}\smallskip}

\newenvironment{lemma}{\smallskip\begin{trivlist}
 \item[\hspace{\labelsep}{\noindent\bf Lemma.}]\it
 }{\end{trivlist}\smallskip}
\newenvironment{propo}{\smallskip\begin{trivlist}
 \item[\hspace{\labelsep}{\noindent\bf Proposition.}]\it
 }{\end{trivlist}\smallskip}
\newenvironment{coro}{\smallskip\begin{trivlist}
 \item[\hspace{\labelsep}{\noindent\bf Corollary.}]\it
 }{\end{trivlist}\smallskip}

\newcommand{\cuadro}{\hfill{$\qed$}}
\newenvironment{examples}{\begin{trivlist}
\item[\hspace{\labelsep}{\noindent\em Examples:}]
}{\end{trivlist}}

\newcommand{\C}{\mathbb{C}}
\newcommand{\D}{\mathbb{D}}

\newcommand{\R}{\mathbb{R}}
\newcommand{\s}{\mathbb{S}}
\newcommand{\T}{\mathbb{T}}
\newcommand{\x}{\mathbb{X}}
\newcommand{\Z}{\mathbb{Z}}
\newcommand{\coh}{\hbox{\rm coh}}
\renewcommand{\deg}{\hbox{\rm deg}}
\newcommand{\degp}{\hbox{\tiny\rm deg}}
\renewcommand{\det}{\hbox{\rm det}}

\newcommand{\id}{\hbox{\rm id}}

\renewcommand{\max}{\hbox{\rm max}}
\renewcommand{\mod}{\hbox{\rm mod}}
\newcommand{\rad}{\hbox{\rm rad}}
\newcommand{\Root}{\hbox{\rm Root}}
\newcommand{\Spec}{\hbox{\rm Spec}}

\def\raya{\raise1.5pt\hbox to 25pt{\vrule height1.5pt depth-1pt
           width25pt}}
\def\rayita{\raise2pt\hbox to 7.5pt{\vrule height1.5pt depth-1pt
           width7.5pt}}
\newcommand{\bulito}{\ {\scriptstyle \bullet}\, \ }
\def\subsetnoteq{\mathbin{\hbox{$\subseteq \joinrel \hskip-8pt \lower3pt
                 \hbox{$\scriptscriptstyle /$}\ $}}}
\def\aderecha#1{\smash{\mathop{\longrightarrow}\limits^{#1}}}

\begin{document}
\title{A CHEBYSHEFF RECURSION FORMULA FOR COXETER POLYNOMIALS}
\author[Helmut Lenzing and J.$\!$~A. de la Pe\~na]{Helmut Lenzing$^1$
and J.$\!$~A. de la Pe\~na$^2$}
\thanks{\quad\\[-12pt]
\noindent $^1$ Institut f\"ur Mathematik, Universit\"at Paderborn.
Germany. helmut@math.uni-paderborn.de}
\thanks{\quad\\[-12pt]
\noindent $^2$ Instituto de Matem\'aticas, UNAM, Circuito
Exterior, Ciudad Universitaria. M\'exico, D.$\!$~F.
jap@matem.unam.mx}

\begin{abstract}
A polynomial $f(T)\in \Z [T]$ is {\em represented\/} by $q(T)\in \Z
[T]$ if $f(T^2)=q^*(T)=T^{\degp\,q}q(T+T^{-1})$; $f(T)$ is {\em
graphically represented\/} if $f(T^2)=\chi^*_M(T)$ for $\chi_M(T)$
the characteristic polynomial of a symmetric matrix $M$. Many
instances of {\em Coxeter polynomials\/} $f_A(T)$, for $A$ a finite
dimensional algebra, are (graphically) representable. We study the
case of extended canonical algebras $A$, see \cite{17}, show that
the corresponding polynomials $f_A(T)$ are representable and satisfy
a Chebysheff type recursion formula. We get consequences for the
eigenvalues of the Coxeter transformation of $A$ showing, for
instance, that at most four eigenvalues may lie outside the unit
circle.
\end{abstract}

\maketitle

\section*{Introduction}
For a finite dimensional $k$-algebra $A$ of finite global dimension,
the {\em Coxeter transformation\/} $\varphi_A$ is an automorphism of
the Grothendieck group $K_0(A)$ such that
$[X^{\bullet}]\varphi_A=[\tau_{D(A)}X^{\bullet}]$ for any complex
$X^{\bullet}$ in the bounded derived category $D(A):=D^b(\mod\,A)$
of finite dimensional left $A$-modules, where $\tau_{D(A)}$ is the
Auslander-Reiten translation in $D(A)$ and $[X^{\bullet}]$ is the
class of $X^{\bullet}$ in $K_0(D(A))\ \aderecha{\sim}\ K_0(A)$. The
characteristic polynomial $f_A(T)=\det\,(T\id -\varphi_A)$, called
the {\em Coxeter polynomial\/} of $A$, and the corresponding
spectrum $\Spec\,\varphi_A=\Root\,f_A(T)$ control the growth
behavior of $\varphi_A$ and hence of $\tau_{D(A)}$, see
\cite{12,11,14}. Clearly, $f_A(T)$ and $\Spec\,\varphi_A$ are
derived invariants of $A$ and establish links between representation
theory and other fundamental areas: Lie algebras, $C^*$-algebras,
spectral theory of graphs among other topics.

A fundamental fact for a hereditary algebra $A=kQ$, when $Q$ is a
{\em bipartite quiver\/} without oriented cycles is that
$\Spec\,\varphi_A\subset \s^1\cup \R^+$. This was shown by A'Campo
\cite{10} as a consequence of the identity $f_A(T^2)=T^n\chi_\Delta
(T+T^{-1})$, where $n$ is the number of vertices of $Q$ and
$\chi_\Delta (T)$ is the characteristic polynomial of the adjacency
matrix of the underlying graph $\Delta$ of $Q$. We shall say that a
polynomial $f(T)\in \Z [T]$ is {\em represented\/} by $q(T)\in \Z
[T]$ if $f(T^2)=q^*(T):=T^{\degp\,q}q(T+T^{-1})$. Moreover, $f(T)$ is
{\em graphically represented\/} by the symmetric matrix $M\in
M_{n\times n}(\Z)$ if $f(T^2)=\chi^*_M(T)$ for the characteristic
polynomial $\chi_M(T)$ of $M$. We shall study representability of
polynomials and obtain consequences for the theory of Coxeter
polynomials.

In the general setting we prove in section~2 that any monic
polynomial $f(T)\in \Z [T]$ with $\Root\,f(T)$ in the unit disk
$\D^1$, with $f(1)\neq0$, is representable.

In section~3 we recall that Coxeter polynomials of several
well-known classes of algebras are representable. We concentrate in
a special new class of algebras, namely the {\em extended canonical
algebras\/} $A$ of the form $A=C[P]$ a one-point extension of a
canonical algebra $C$ by an indecomposable projective $C$-module
$P$. This class of algebras was introduced and studied in \cite{17}
(background about canonical algebras and their $K$-theory may be
seen in \cite{12,15,9,20}). For a canonical algebra $C$ of type
$(p_1,p_2,\ldots,p_t)$ (here $t\ge 2$, $p_1,\ldots,p_t$ are integers
$\ge 1$) we associate a unique Coxeter polynomial
$\hat{f}_{(p_1,\ldots,p_t)}(T)$ corresponding to an extended
canonical algebra $A=C[P]$. We shall prove that:

\medskip

$\bulito$ $\hat{f}_{(p_1,\ldots,p_t)}(T^2)=q^*_{(p_1,\ldots,p_t)}(T)$
for a polynomial $q_{(p_1,\ldots,p_t)}(T)\in \Z [T]$

$\bulito$ $q_{(p_1,\ldots,p_t+1)}(T)=Tq_{(p_1,\ldots,p_t)}(T)
-q_{(p_1,\ldots,p_t-1)}(T)$.

\medskip

The above recursion formula has the shape of the {\em Chebysheff
formula\/} defining the (normalized) Chebysheff polynomials of the
second kind. As consequence of the recursion formula and Sturm's
theorem we are able to prove essential properties of
$\Spec\,\varphi_A$ which are used in \cite{17}. For instance, we show
that $\varphi_A$ accepts at most $4$ eigenvalues outside $\s^1$.

We believe that representability of Coxeter polynomials is an
important feature with the potential for further applications of
spectral theory to representation theory of algebras.

The research for this paper was done during a visit of the second
named author to Paderborn.

\section{Self-reciprocal and representable polynomials} 

\subsection{} 
Let $p(T)=a_0T^n+a_1T^{n-1}+\cdots +a_{n-1}T+a_n$ be an integral
polynomial of degree $\deg\,p=n$. We write $\varepsilon (p)\in
\{0,1\}$ to indicate the {\em parity\/} of $\deg\,p$, that is,
$\varepsilon (p)\equiv \deg\,p\,(\mod\,2)$. We shall deal mainly with
{\em monic\/} polynomials, where $a_0=1$ holds.

We recall that $p(T)$ is said to be {\em self-reciprocal\/} if
$a_i=a_{n-1}$ for $i=0,\ldots,n$. For any polynomial $q(T)$, the
polynomial $q^*(T)=T^{\degp\,q}q(T+T^{-1})$ is self-reciprocal, it
will be called the {\em symmetrization polynomial\/} of $q(T)$. The
following is elementary:

\begin{lemma} The following are equivalent for the polynomial
$p(T)\in \Z [T]$:

\begin{itemize}
\item[{\rm (a)}] $p(T)$ is self-reciprocal

\item[{\rm (b)}] $p(T)=T^{\degp\,p}p\left(\frac{1}{T}\right)$

\item[{\rm (c)}] $p(T)=(T-1)^{2s}q(T)$ for some $s\ge 0$ and some
polynomial $q(T)\in \Z [T]$ such that $q(1)\ne 0$ and for any $0\ne
\lambda \in \Root\,q(T)$, then $\lambda^{-1}\in \Root\,q(T)$.\cuadro
\end{itemize}
\end{lemma}

\medskip

If $p(T)$ is self-reciprocal observe the following:

\medskip

(i) if $p(T)$ has odd degree, then $T+1$ is a divisor of $p(T)$.

(ii) if $p(T)$ has even degree $2k$, we find a polynomial $q(T)\in \Z
[T]$ such that $p(T)=q^*(T)$ (in fact, observe that
$q(T)=b_0T^k+b_1T^{k-1}+\cdots +b_k$ and if we have fixed
$b_0=a_0,b_1,\ldots,b_i$ with $i<k$, then $a_i=b_i+\sum\limits_{0\le
i-2j}b_{i-2j}\left({n-i+2j\atop j}\right)$; hence $q(T)\in \Z [T]$).

\begin{coro}
For every self-reciprocal polynomial $p(T)\in \Z [T]$, there is a
polynomial $q(T)\in \Z [T]$, with $2\deg\,q+\varepsilon (p)=\deg\,p$,
such that $p(T)=(T+1)^{\varepsilon (p)}q^*(T)$.\cuadro
\end{coro}

\subsection{} 
By (1.1) we know that representable polynomials are self-reciprocal.
The converse partially holds as follows from the next Proposition.

For $0\ne f(T)\in \Z [T^2]$ consider $\hat{f}(T)\in \Z [T]$ a
polynomial such that $\hat{f}(T^2)=T^{\degp\,f}f(T+T^{-1})(=f^*(T))$.

\begin{propo}
The map $f\mapsto \hat{f}$ establishes an isomorphism of
multiplicative semigroups between $\Z [T^2]\setminus \{0\}$ and the
non-zero self-reciprocal polynomials in $\Z [T^2]$.
\end{propo}

\begin{proof}
Observe that the map is well-defined, it is injective and preserves
multiplication. Let $p(T)$ be a self-reciprocal polynomial in $\Z
[T^2]$, then $p(T)=q^*(T)$ for some $q(T)\in \Z [T]$, by (1.1). Hence
$f(T)=q(T^2)\in \Z [T^2]$ with $\hat{f}(T^2)=f^*(T)=q^*(T^2)=p(T^2)$
and $\hat{f}(T)=p(T)$.
\end{proof}

\subsection{} 
We consider the class of {\em cyclotomic polynomials\/} $\phi_n(T)\in
\Z [T]$. We recall that $\phi_1(T)=T-1$ and for $n>1$,
$$\phi_n(T)\prod_{1<d<n\atop d|n}\phi_d(T)=T^{n-1}+\cdots
+T^2+T+1=v_{n}(T)=\frac{T^n-1}{T-1}.$$ The roots of $\phi_n(T)$ are
the $n$-th primitive roots of $1$ in $\C$, hence the polynomial
$\phi_n(T)$, for $n\ge 2$ is self-reciprocal (since with $\lambda
\in \Root\,\phi_n(T)$, also $\lambda^{-1}$ is a $n$-th primitive
root of $1$). Only $\phi_1(T)=T-1$ and $\phi_2(T)=T+1$ have odd
degree among the $\phi_n(T)$.

The following simple remark is useful. Consider the multiplicative
semigroup $A$ of self-reciprocal polynomials. The irreducible
elements $p(T)\in A$ will be called {\em $sr$-irreducible\/}
polynomials. Clearly, there are two types of $sr$-irreducible
polynomials:

\begin{itemize}
\item[{\rm (i)}] $p(T)$ irreducible in $\Z [T]$ and $p(T)\in A$;

\item[{\rm (ii)}] $p(T)\bar{p}(T)$, where $p(T)$ is the minimal
polynomial with $p(\lambda)=0$ for certain $\lambda \in \C$ such that
$p(\lambda^{-1})\ne 0$ and $\bar{p}(T)$ is the minimal polynomial
with $\bar{p}(\lambda^{-1})=0$.
\end{itemize}

\begin{lemma}
The following holds:

\begin{itemize}
\item[{\rm (a)}] There are unique decompositions into $sr$-irreducible
polynomials in $A$.

\item[{\rm (b)}] $\phi_n(T^2)$ is an irreducible polynomials in $\Z
[T]$, for $n\ge 2$.
\end{itemize}
\end{lemma}

\begin{proof}
(a) follows from the remarks above.

(b): That $\phi_n(T)$ is irreducible is well-known \cite{1}. Assume
$\phi_n(T^2)=g_1(T)g_2(T)$ for two polynomials $1\ne g_1(T),g_2(T)\in
\Z [T]$. Assume $g_1(T)=T^s+a_1T^{s-1}+\cdots +a_s$ and
$g_2(T)=T^t+\cdots +b_t$. If $1$ or $-1$ are roots of $g_i(T)$, then
also $1$ would be a root of $\phi_n(T)$. In that case $n=1$. We may
assume that $g_j(T)$ (as product of cyclotomic polynomials) has even
degree, $j=1,2$. Assume $1\le i\le s$ is a maximal odd index with
$a_i\ne 0$, then $a_i$ is the coefficient of $T^{t+i}$ in
$g_1(T)g_2(T)$, a contradiction. Therefore $g_i(T)=h_i(T^2)$ for some
integral polynomial $h_i(T)$, $i=1,2$. Then $\phi_n(T)=h_1(T)h_2(T)$
yields the result.
\end{proof}

\subsection{} 
Our interest in the representability of polynomials is due to the
important interrelations of $\Root\,p(T)$ and $\Root\,q(T)$ when
$p(T^2)=q^*(T)$. A first glimpse on these relations is the following
summary of known facts.

\begin{propo}
Let $p(T)$ and $q(T)\in \Z [T]$ be such that $p(T^2)=q^*(T)$. Then

\begin{itemize}
\item[{\rm (a)}] $\Root\,q(T)$ is symmetric with respect to $0$, that
is, if $\lambda \in \Root\,q(T)$, also $-\lambda \in \Root\,q(T)$.

\item[{\rm (b)}] $\mu \in \Root\,p(T)$, then
$\mu^{-1},\bar{\mu},\bar{\mu}^{-1}\in \Root\,p(T)$.

\item[{\rm (c)}] $\Root\,p(T)\subset \s^1\cup \R^+$ (where
$\s^1=\{v\in \C\colon \Vert v\Vert =1\}$ is the unit circle and
$\R^+$ the positive real numbers, $0\in \R^+$) if and only if
$\Root\,q(T)\subset \R$.

\item[{\rm (d)}] $\Root\,p(T)\subset \s^1$ (resp. $\s^1\setminus
\{1\}$) if and only if $\Root\,q(T)\subset [-2,2]$ (resp. $(-2,2)$).
\end{itemize}
\end{propo}

\begin{proof}
Observe that $0\ne \lambda\in \C$, then $\mu=\lambda^2\in\Root\,p(T)$
if and only if $\lambda +\lambda^{-1}\in\Root\,q(T)$. Assume $\lambda
=r(a+ib)$ with $r\in \R^+$, $a^2+b^2=1$, then $\lambda +
\lambda^{-1}=(r+r^{-1})a+i(r-r^{-1})b$.

(a): Clear; (b): if $\mu =\lambda^2\in\Root\,p(T)$, then $\mu^{-1}$
yields $\lambda +\lambda^{-1}\in \Root\,q(T)$, then $\mu^{-1}\in
\Root\,p(T)$. Since $p(T)$ has real coefficients, then $\bar{\mu}$
(and therefore $\bar{\mu}^{-1}$) is in $\Root\,p(T)$.

(c): Observe $\mu \in \Root\,(T)$ lies in $\s^1$ (resp. in $\R^+$)
if and only if $r=1$ (resp. $b=0$) if and only if $\lambda
+\lambda^{-1}\in \R$.

(d): In the above case $\lambda +\lambda^{-1}\in [-2,2]$ if and only
if $b=0$ if and only if $\mu \in \s^1$.
\end{proof}

\subsection{} 
Certain type of polynomials are of special interest for us. We say
that a polynomial $q(T)\in \Z [T]$ is of {\em graphical type\/} if
$q(T)=\chi_M(T)$ for a $n\times n$ integral symmetric matrix $M$,
where $\chi_M(T)$ is the {\em characteristic polynomial\/} of $M$,
that is,
$$\chi_M(T)=\det\,(TI_n-M).$$

We say that a polynomial $p(T)\in \Z [T]$ is {\em graphically
represented\/} if $p(T^2)=\chi^*_M(T)$ for a symmetric matrix $M\in
M_{n\times n}(\Z)$.

\begin{examples}
(a) The polynomial $q(T)=T^2-1$ is not of graphical type. In
particular, $\phi_4(T)$ is not graphically represented.

(b) For a finite graph $\Delta$ we associate the {\em adjacency
matrix\/} $A(\Delta)=(a_{ij})$ of size $n\times n$ with
$\Delta_0=\{1,\ldots,n\}$ the vertices of $\Delta$ and $a_{ij}$ the
number of edges between $i$ and $j$. Of particular interest are the
following graphs:

\medskip

\renewcommand{\AA}{\mathbb{A}}
\newcommand{\EE}{\mathbb{E}}
\newcommand{\DD}{\mathbb{D}}
\newcommand{\TT}{\mathbb{T}}

\newcommand{\dynkin}{$\xymatrix@R2pt@C10pt{
 {\AA_n}:& 1\ar@{-}[r]&2\ar@{-}[r]&\cdots\ar@{-}[r]&n && {\EE_p}:&1\ar@{-}[r]&2\ar@{-}[r] &3\ar@{-}[r]&5\ar@{-}[r]&\cdots\ar@{-}[r]&p &(p=6,\,7,\,8).\\
         &1\ar@{-}[dr]&           &                &  &&         &           &            &4\ar@{-}[u]&     &         &&\\
 {\DD_n}:&            &3\ar@{-}[r]&\cdots\ar@{-}[r]&n &&          &          &            &           &                   &         &&\\
         &2\ar@{-}[ur]&           &                &  &&           &         &            &           &                   &         &&\\
}$}

\newcommand{\kron}{$\xymatrix@R10pt@C10pt{ 1\ar@/ ^1pc/[rr]\ar@/ _1pc/[rr]&{\vdots}&2}$}

\newcommand{\stars}{$\xymatrix@R2pt@C15pt{
                        &                        &(1,2)\ar@{-}[r]            &(1,3)\ar@{-}[r] &\cdots\ar@{-}[r]&(1,p_1)\\
                        &                        &(2,2)\ar@{-}[r]            &(2,3)\ar@{-}[r] &\cdots\ar@{-}[r]&(2,p_2)        \\
 \TT_{(p_1,\ldots,p_t)}: &1\ar@{-}[ruu]\ar@{-}[ru]&\vdots                     &                &                & \vdots  \\
                        &                        &(t,2)\ar@{-}[lu]\ar@{-}[r] &(t,3)\ar@{-}[r]&\cdots\ar@{-}[r]&(t,p_t)\\
}$}

$\bulito$ {\em Dynkin\/} graphs:
\begin{figure}[h]
\begin{center}
\dynkin
\end{center}
\end{figure}

$\bulito$ the graph $K_s\colon
\xymatrix{1\ar@/^/@{-}[r]&2\ar@/^/@{-}[l]}
\hskip-26pt\lower1.5pt\hbox{\vdots}$
\hskip26pt with $s$ edges

\medskip

$\bulito$ the {\em stars:}
\begin{figure}[h]
\stars
\end{figure}
\end{examples}
which will also be denoted $[p_1,\ldots,p_t]$.

For simplicity we write $\chi_{A(\Delta)}=:\chi_\Delta$.

For consideration on the spectra of $A(\Delta)$
($=\Root\,\chi_{A(\Delta)}$) we refer to \cite{5,6}. In particular we
mention the following important facts:

(i) If $M$ is a symmetric $n\times n$-matrix, then $\chi_M(T)$ is a
self-reciprocal polynomial and $\Root\,\chi_M(T)\subset \R$.

(ii) $M=A(\Delta)$ is an irreducible matrix if and only if $\Delta$
is connected. In that case, there is a root $\rho (\Delta)$ of
$\chi_\Delta (T)$ called the {\em spectral radius\/} of $M$ (and of
$\Delta$) such that
$$\rho (\Delta)=\max\,\{\Vert \lambda \Vert \colon \lambda \in
\Root\,\chi_\Delta\}.$$
Moreover, $\rho (\Delta)$ is a simple root of $\chi_\Delta (T)$.

The heart of the argument of A'Campo mentioned in the Introduction is
the following:

\begin{propo}
If $p(T)$ is graphically represented by $q(T)$, then
$\Root\,q(T)\subset \R$ and $\Root\,p(T)\subset \s^1\cup \R^+$.
\end{propo}

\begin{proof}
Assume $q(T)=\chi_M(T)$ for a symmetric matrix $M$ and
$p(T^2)=q^*(T)$. Then $\Root\,q(T)\subset \R$ and by (1.4),
$\Root\,p(T)\subset \s^1\cup \R^+$.
\end{proof}

\section{Integral polynomials with roots in $\s^1$} 

\subsection{} 
We show representability of the cyclotomic polynomials $\phi_n(T)$,
that is, we want to identify those $f_n(T)\in \Z [T]$ such that
$f^*_n(T)=\phi_n(T^2)$. For that purpose recall (see for example
\cite{1,2}) that the (normalized) {\em Chebysheff polynomials\/} (of
the second kind) $u_n(T)\in \Z [T]$ are given by
$$u_n(T)=2U_n(T/2)\hbox{ where
}U_n(T)=\frac{\sin\,n\theta}{\sin\,\theta}\hbox{ and }T
=\cos\,\theta .$$
These polynomials satisfy $\hat{u}_n(T)=v_n(T)$. Moreover the family
$(u_n(T))_n$ may be inductively constructed by the rules:
$$u_0(T)=1,\ u_1(T)=T\hbox{ and}$$
$$u_{n+1}(T)=Tu_n(T)-u_{n-1}(T)\hbox{ for }n\ge 1.$$

\subsection{} 
The characteristic polynomial $\chi_{[p_1,\ldots,p_t]}$ of stars
$\mathbb{T}_[p_1,\ldots,p_t]$ or other trees $\Delta$ may be
constructed by repeatedly applying the following device: if $a$ is a
vertex in $\Delta$ with a unique neighbor $b$ and $\Delta'$ (resp.\
$\Delta''$) is the full subgraph of $\Delta$ with vertices
$\Delta_0\setminus\{a\}$ (resp.\ $\Delta_0\setminus\{a,b\}$), then
$$
\chi_\Delta(T)=T\,\chi_{\Delta'}(T)-\chi_{\Delta''}(T).
$$
The following is a well known consequence.
\begin{lemma}
The characteristic polynomial of the linear graph $\mathbb{A}_n=[n]$
is the $n$-th normalized Chebysheff polynomial $u_n$. Moreover,
$v_{n+1}(T^2)=u_n^*(T)$.\cuadro
\end{lemma}

\subsection{}Observe that $\phi_1(T)=T-1$ is not representable.
\begin{propo}
For each $n\ge 2$, there is an irreducible factor $f(T)$ of
$u_{n-1}(T)$ such that $\phi_n(T^2)=f^*(T)$.
\end{propo}

\begin{proof}
For $n=2$, we have $\phi_2(T^2)=T^2+1=u^*_1(T)$.

Assume $f^*_i(T)=\phi_i(T^2)$ for $i<n$, then
$$\phi_n(T^2)\prod_{1<d<n\atop d|n}\phi_d(T^2)=v_{n}(T^2)=
u^*_{n-1}(T).$$
Using (1.4), consider a decomposition
$u_{n-1}(T)=\prod\limits^s_{i=1}g_i(T)$ in $\Z [T]$ such that $g_i(T)$
is $sr$-irreducible. By induction hypothesis we may assume that
$\phi_d(T^2)=f^*_d(T)$ for certain $f_d(T)=g_i(T)$, $1\le i\le s$.
The result for $\phi_n$ follows from (1.4).
\end{proof}

As concrete examples we calculate:

$\phi_3(T)=T^2+T+1$ with $\phi_3(T^2)=u^*_2(T)$; $\phi_4(T)=T^2+1$
with $\phi_4(T^2)=f^*_4(T)$ for $f_4(T)=T^2-2$ which is a factor of
$u_3(T)=T^3-2T$; $\phi_5(T)=T^4+T^3+T^2+T+1$ with
$\phi_5(T^2)=u^*_4(T)$.

\subsection{} 
A famous result by Kronecker \cite{3} (see also \cite{4}) states that
any monic polynomial $p(T)\in \Z [T]$ with $\Root\,p(T)\subset
\D^1=\{v\in \C \colon \Vert v\Vert \le 1\}$ is a product of
cyclotomic polynomials (and hence all roots are roots of unity). As
consequence of the above considerations, we get the following

\begin{theor}
Let $p(T)\in \Z [T]$ be a monic polynomial with $\Root\,p(T)\subset
\D^1\setminus\{1\}$. Then there exists a polynomial $q(T)\in \Z [T]$
such that $p(T^2)=q^*(T)$.\cuadro
\end{theor}

\section{Representability and the Coxeter polynomial} 

\subsection{} 
Let $A$ be a finite dimensional $k$-algebra with $k$ an algebraically
closed field. For simplicity we assume $A=kQ/I$ for a quiver without
oriented cycles and $I$ an ideal of the path algebra. We refer to
\cite{9,10,11,12} for concepts and examples.

We recall that the {\em Coxeter transformation\/} $\varphi_A$ of the
Grothendieck group $K_0(A)$ is induced by the Auslander-Reiten
translation $\tau_A$ of the bounded derived category $D^b(\mod\,A)$
of the finite dimensional $A$-modules $\mod\,A$. This $n\times
n$-integral matrix plays an important role in the study of the
representation theory of certain algebras. By $f_A(T)$ we denote the
characteristic polynomial of $\varphi_A$, called the {\em Coxeter
polynomial\/} of $A$.

We substantiate the consideration of $\varphi_A$:

(i) Let $S_1,\ldots,S_n$ be a complete system of pairwise
non-isomorphic simple $A$-modules, $P_1,\ldots,P_n$ the corresponding
projective covers and $I_1,\ldots,I_n$ the injective envelopes. Then
$\varphi_A$ is defined by $[P_i]\varphi_A=-[I_i]$, where $[X]$
denotes the class of a module $X$ in $K_0(A)$.

(ii) For a hereditary algebra $A=kQ$, the {\em spectral radius\/}
$\rho (\varphi_A)$ determines the representation type of $A$ in the
following manner:

\begin{itemize}
\item[$\bulito$] $A$ is representation-finite if $1=\rho (\varphi_A)$
is not a root of $f_A(T)$.

\item[$\bulito$] $A$ is tame if $1=\rho
(\varphi_A)\in\Root\,\chi_{\varphi_A}(T)$.

\item[$\bulito$] $A$ is wild if $1<\rho (\varphi_A)$. Moreover, in
this case, Ringel \cite{13} shows that $\rho (A)$ is a simple root of
$f_A(T)$. For a non-preprojective indecomposable module $X$ the
vectors $[\tau^n_AX]$ grow exponentially with $\rho (\varphi_A)$
(that is, $\lim\limits_{n\to \infty}\frac{\|[\tau^n_AX]\|}{\rho
(\varphi_A)^n}>0$), see \cite{14}.
\end{itemize}

(iii) In the particular case $A=kQ$ and $Q$ is a bipartite quiver, as
we recalled in the introduction, $f_A(T)$ is graphically represented
by the underlying graph of $Q$, that is:
$$f_A(T^2)=\chi^*_\Delta (T)$$
where $\Delta =|Q|$.

(iv) In \cite{17} the authors introduced {\em supercanonical
algebras\/} as generalization of other important classes of algebras.
Let $S_1,\ldots,S_n$ be a finite number of posets, $t\ge 2$ and
numbers $\lambda_3,\ldots,\lambda_t\in k\setminus \{0\}$ pairwise
different. The supercanonical algebra
$A=A(S_1,\ldots,S_t;\lambda_3,\ldots,\lambda_t)$ is the quiver
algebra with a unique source $\alpha$, a unique sink $\omega$ as in
the picture:

\newcommand{\al}{\alpha}
\newcommand{\om}{\omega}

\newcommand{\supercan}{$\xymatrix{
                  & \al\ar[dl]\ar[d]\ar[dr]\ar[drr]&               &         \\
 *+[F]{S_1}\ar[dr]& *+[F]{S_2}\ar[d]               &\cdots\ar[dl] &*+[F]{S_t}\ar[dll]\\
                  & \om                            &              &
} $}

\begin{figure}[h]
\begin{center}
\supercan
\end{center}
\end{figure}
\eject %

\noindent all paths from $\alpha$ to $\omega$ passing through $S_i$
are equal in $A$ to a unique path $\kappa_i$, additionally in $A$
the $t-2$ relations $\kappa_i=\kappa_2-\lambda_i\kappa_1$,
$i=3,\ldots,t$ are satisfied.

In case, $S_1,\ldots,S_t$ are linear quivers
$$S_i=[p_i-1]\colon 1\to 2\to \cdots \to p_i-1$$
the algebra $A(S_1,\ldots,S_t;\lambda_1,\ldots,\lambda_t)$ denoted by
$C(p_1,\ldots,p_t;\lambda_3,\ldots,\lambda_t)$ is called a {\em
canonical algebra}. See \cite{9,11,15}.

Since $A=A(S_1,\ldots,S_t;\lambda_1,\ldots,\lambda_t)$ is a one-point
extension of the poset algebra $B=A/(\alpha)$, the Coxeter
transformations are related by formulas:
$$\varphi_A=\left[\begin{array}{ccc}
\varphi_B &| &-C^{-t}_Bv^t\\
\hbox{---}\quad \hbox{---} &\hbox{---} &\hbox{---}\quad \hbox{---} \\
-v\varphi_B &| &q_B(v)-1
\end{array}\right]$$
where $v=[M]\in K_0(B)$ is such that $\rad\,P_\alpha =M$; $C_B$ is
the Cartan matrix satisfying $\varphi_B=-C^{-t}_BC_B$ and $q_B\colon
K_0(B)\to \Z$ is the quadratic form defined as
$q_B(x)=x(C^{-1}_B+C^{-t}_B)x^t$. A simple calculation shows that
$$f_A(T)=(T-1)^2\prod^t_{i=1}f_{S_i}(T)$$
where $\varphi_{S_i}$ is the Coxeter transformation of the poset
$S_i$, $i=1,\ldots,t$.

(iv) A tool to deal with the calculation of Coxeter
transformations is the following {\em reduction formula\/}
\cite[18.3.3]{12} (see also \cite{18} for generalizations and special
cases):

If $A=B[P]$ is a one-point extension of an algebra $B$ with an
indecomposable projective module $P$ associated to a source $b$ in
$B$. Let $C=B/(b)$, then
$$f_A(T)=(1+T)f_B(T)-Tf_C(T).$$
Of particular interest is the hereditary case where any algebra can
be constructed by repeated one-point extensions using only projective
modules. For the Coxeter polynomial of the star
$H=\T_{p_1,\ldots,p_t}$ the formula yields:
$$f_{[p_1,\ldots,p_t]}(T):=f_H(T)=
(T+1)\prod^t_{i=1}v_{p_i}(T)-T\sum\limits^t_{i=1}v_{p_i-1}(T)
\prod_{i\ne j}v_{p_j}(T)$$ where $v_n=(T^n-1)/(T-1)$ is the Coxeter
polynomial of the linear quiver $[n-1]$.

\subsection{} 
Let $C=C(p_1,\ldots,p_t;\lambda_3,\ldots,\lambda_t)$ be a canonical
algebra. In \cite{17} the authors proposed to study an interesting
class of algebras according with the following result:

\begin{propo}
The derived equivalence class of the one-point extensions of $C$ by
an indecomposable projective or injective module is independent of
the particular choice of the module.
\end{propo}

\begin{proof}
If $P_1$ and $P_2$ are indecomposable projective modules over $C$,
there is an automorphism of $D^b(\mod\,C)=D^b(\coh\,\x)$ sending
$P_1$ to $P_2$, see \cite{20}. Here $\coh\,\x$ denotes the coherent
sheaves over a weighted projective line $\x$. The assertion now
follows from \cite{20}.
\end{proof}

We call an algebra of the form $C[P]$ with $P$ indecomposable
projective, and {\em extended canonical algebra\/} of type
$(p_1,\ldots,p_t;\lambda_3,\ldots,\lambda_t)$. The Coxeter
transformation $f_{\varphi_{\lower4pt\hbox{$\scriptscriptstyle
C[P]$}}}(T)=:\hat{f}_{(p_1,\ldots,p_t)}(T)$ only depends on the
numbers $p_1,\ldots,p_t$.

\begin{coro}
The extended canonical algebra, as defined above, has Coxeter
polynomial
$$\hat{f}_{(p_1,\ldots,p_t)}(T)=(T+1)(T-1)^2\prod^t_{i=1}v_{p_i}(T)
-Tf_{[p_1,\ldots,p_t]}(T).$$
\vskip-28pt \cuadro
\end{coro}

\subsection{} 
There are interesting phenomena happening for extended canonical
algebras and their Coxeter transformation. A couple of {\em
examples\/} taken from \cite{17}:

(a) Let $C$ be a canonical tubular algebra (i.e. the star
$\T_{(p_1,\ldots,p_t)}$ is of extended Dynkin type). Then the
extended canonical algebra $C[P]$ is of canonical derived type
$(p_1,\ldots,p_{t-1},p_t+1)$.

(b) There is a {\em finite list of types\/} which is critical
(according to the lexicographical order) for extended canonical
algebras to have the spectral radius of their Coxeter
transformations equal to $1$. As a consequence there are only
finitely many types $(p_1,\ldots,p_t)$ such that the corresponding
extended canonical algebra $A$ has $\rho (\varphi_A)=1$. (A result
which is proved as consequence of the considerations in this paper,
namely Theorem~3.6).

(c) With two exceptions, namely the types $(3,3,3,3)$ and
$(2,2,2,2,4)$, all extended canonical algebras $A$ with
$\rho(\varphi_A)=1$ have a periodic Coxeter transformation
$\varphi_A$.

We propose to use the theory of representability of polynomials as
developed before in order to study the eigenvalues of Coxeter
matrices of extended canonical algebras.

\subsection{\it Fix notation:} 
Let $H$ be the star of type $[p_1,\ldots,p_t]$ with Coxeter
polynomial $f_{[p_1,\ldots,p_t]}(T)$; let $C$ be the canonical
algebra of type $(p_1,\ldots,p_t;\lambda_1,\ldots,\lambda_t)$ with
Coxeter polynomial $f_{(p_1,\ldots,p_t)}(T)$; let $A=C[P]$ be the
extended canonical algebra with Coxeter polynomial
$\hat{f}_{(p_1,\ldots,p_t)}(T)$.

\begin{propo}
The above Coxeter polynomials are representable in the following way:

\begin{itemize}
\item[{\rm (a)}] $f_{[p_1,\ldots,p_t]}(T^2)=
\chi^*_{[p_1,\ldots,p_t]}(T)$, where $\chi_{[p_1,\ldots,p_t]}(T)$
denotes the characteristic polynomial of the star
$\T_{(p_1,\ldots,p_t)}$

\item[{\rm (b)}] $f_{(p_1,\ldots,p_t)}(T^2)=
\chi^*_{K_2}(T)\prod\limits^t_{i=1}\chi^*_{[p_i-1]}(T)$, where
$\chi_{[p]}(T)$ denotes the characteristic polynomial of the linear
path $[p]$ and $K_2$ is the Kronecker diagram $\xymatrix{\bulito
\ar@{-}@/^/[r]&\bulito \ar@{-}@/^/[l]}$

\item[{\rm (c)}] $\hat{f}_{(p_1,\ldots,p_t)}(T^2)=
\left(T\chi_{K_2}(T)\prod^t_{i=1}\chi_{[p_i-1]}(T)-
\chi_{[p_1,\ldots,p_t]}(T)\right)^*$.
\end{itemize}

Moreover for $p_t\ge 2$ the following formula holds:

\begin{itemize}
\item[{\rm (d)}]$\hat{f}_{(p_1,\ldots,p_{t+1})}(T)=
T\hat{f}_{(p_1,\ldots,p_t)}(T)-\hat{f}_{(p_1,\ldots,p_t-1)}(T)$,
where $\hat{f}_{(p_1,\ldots,p_{t-1},0)}(T)$ is to be understood as
$\hat{f}_{(p_1,\ldots,p_{t-1}-1)}(T)$.
\end{itemize}
\end{propo}

\begin{proof}
(a): Since the star $\T_{p_1,\ldots,p_t}$ is a bipartite quiver, the
formula follows from (1.5).

(b): Clearly, $\chi^*_{K_2}(T)=T^4-2T^2+1=(T^2-1)^2$ and
$f_{[p_i]}(T^2)=\chi^*_{[p_i]}(T)$ by (1.5).

(c): By Corollary (3.2),
$$\hat{f}_{(p_1,\ldots,p_t)}(T)=(T+1)f_{(p_1,\ldots,p_t)}(T)-
Tf_{[p_1,\ldots,p_t]}(T),$$
then
\begin{eqnarray*}
\lefteqn{\hat{f}_{(p_1,\ldots,p_t)}(T^2)=}\\
&=&(T^2+1)\chi^*_{K_2}(T)\prod^t_{i=1}\chi^*_{[p_i-1]}(T)
-T^2\chi^*_{[p_1,\ldots,p_t]}(T)=\\
&=&T^n\left[(T+T^{-1})\chi_{K_2}(T+T^{-1})
\prod^t_{i=1}\chi_{[p_i-1]}(T+T^{-1})
-\chi_{[p_1,\ldots,p_t]}(T+T^{-1})\right]
\end{eqnarray*}

(d): Since $\chi_{[p_1,\ldots,p_{t+1}]}(T)=T\chi_{[p_1,\ldots,p_t]}(T)
-\chi_{[p_1,\ldots,p_{t-1}]}(T)$ by (1.6.v), and also
$$\chi_{[p_t]}(T)=T\chi_{[p_t-1]}(t)-\chi_{[p_t-2]}(T),$$
we get from (c) formula (d).
\end{proof}

\subsection{} 
Consider the integral polynomial
$$q_{(p_1,\ldots,p_t)}(T)=T\chi_{K_2}(T)\prod^t_{i=1}\chi_{[p_i-1]}(T)
-\chi_{[p_1,\ldots,p_t]}(T)$$
satisfying
$$\hat{f}_{(p_1,\ldots,p_t)}(T^2)=q^*_{(p_1,\ldots,p_t)}(T).$$

\begin{lemma}
$\phantom{a}$

\begin{itemize}
\item[{\rm (a)}] For $p_t\ge 2$, the following holds:
$$q_{(p_1,\ldots,p_t+1)}(T)=
Tq_{(p_1,\ldots,p_t)}(T)-q_{(p_1,\ldots,p_t-1)}(T)$$

\item[{\rm (b)}] The polynomials $q_{(p_1,\ldots,p_t)}(T)$ satisfy
the following conditions (we consider the lexicographical order of
the indices):

\begin{itemize}
\item[{\rm $\bulito$}] $q_{(p_1,\ldots,p_t)}$ is a monic polynomial
of degree $\deg\,q_{(p_1,\ldots,p_t)}=1+\deg\,q_{(p_1,\ldots,p_t-1)}$;

\item[{\rm $\bulito$}] any two successive functions among
$q_{(p_1,\ldots,p_t-1)}(T)$, $q_{(p_1,\ldots,p_t)}(T)$,\break
$q_{(p_1,\ldots,p_t+1)}(T)$ have no common root;

\item[{\rm $\bulito$}] if one of the functions vanishes in $\lambda
\in \R$, then the immediate successor and predecessor of that
function take real values of different sign at $\lambda$ (that is,
if $q_{(p_1,\ldots,p_t)}(\lambda)=0$, then both
$q_{(p_1,\ldots,p_t+1)}(\lambda)$ and
$q_{(p_1,\ldots,p_t-1)}(\lambda)$ are real and
$q_{(p_1,\ldots,p_t+1)}(\lambda)q_{(p_1,\ldots,p_t-1)}(\lambda)<0$).
\end{itemize}
\end{itemize}
\end{lemma}

\begin{proof}
(b): A common zero $\lambda$ of two functions in the series is a
common zero of the whole series. Hence $\lambda$ is a zero of
$$q_{(1,\ldots,1)}(T)=T\chi_{K_2}(T)-\chi_{[1,\ldots,1]}(T)=
T(T^2-4)-T$$
Then either $\lambda =0$ or $\pm \sqrt{5}$. But then $\lambda$ is
also root of
$$\chi_{[1,\ldots,1,2]}(T)=T^2-1,$$
a contradiction.

If $q_{(p_1,\ldots,p_t)}(\lambda)=0$ for some $\lambda \in \R$. Then
as before $\lambda \ne 0$ (else $\chi_{[1,\ldots,1,2]}(0)=0$), hence
$q_{(p_1,\ldots,p_t+1)}(\lambda)=-q_{(p_1,\ldots,p_t-1)}(\lambda)$.
\end{proof}

Under the conditions shown in the lemma, a version of {\em Sturm's
Theorem\/} \cite[\S~20]{19} assures that given any interval $[\alpha
,\beta]\subset \R$ and the roots $\lambda_1\le \cdots \le \lambda_s$
of $q_{(p_1,\ldots,p_{t+1})}(T)$ in $[\alpha,\beta]$, then
$q_{(p_1,\ldots,p_t)}(T)$ has roots $\lambda'_1\le \cdots \le
\lambda'_{s-1}$ in $[\alpha,\beta]$ satisfying $\lambda_1\le
\lambda'_1\le \lambda_2\le \lambda'_2\le \cdots \le \lambda_{s-1}\le
\lambda'_{s-1}\le \lambda_s$ (interlacing property). Moreover, all the
real roots of $q_{(p_1,\ldots,p_t)}(T)$ are simple roots.

\subsection{Theorem.} 
{\it Consider the extended canonical algebra $A$ of type
$(p_1,\ldots,p_{t-1},p_t+1)$ and $A'$ an extended canonical algebra
of type $(p_1,\ldots,p_t)$. Then the following holds:

\begin{itemize}
\item[{\rm (a)}] If $\rho (\varphi_A)=1$, then also $\rho
(\varphi_{A'})=1$.

\item[{\rm (b)}] $\varphi_A$ accepts at most $4$ eigenvalues outside
$\s^1$.
\end{itemize}
}

\begin{proof}
(a): If $\rho (\varphi_A)=1$, then $\chi_{\varphi_A}(T)$ has all its
roots in $\s^1$, then by (1.4)\break
$\Root\,q_{(p_1,\ldots,p_{t+1})}(T)\subset [-2,2]$. By Sturm's
Theorem, $\Root\,q_{(p_1,\ldots,p_t)}(T)\subset [-2,2]$ and hence
$\Root\,\chi_{\varphi_{A'}}(T)\subset \s^1$.

(b): By induction on the sum of the weights $\sum\limits^t_{i=1}p_i$,
we may assume that the Coxeter transformation $\varphi_{A'}$ of the
extended canonical algebra of type $(p_1,\ldots,p_t)$ accepts at most
$4$ roots outside $\s^1$. We shall prove that $\varphi_A$ for $A$ of
type $(p_1,\ldots,p_t+1)$ has the same property. According to (1.7)
and (3.5), the roots of $q_{(p_1,\ldots,p_t)}(T)$ can be enumerated
as $-2\le \lambda_1<\lambda_2<\cdots <\lambda_{s-1}<\lambda_s\le 2$
with $s\ge n-4$ where $n=3-t+\sum\limits^t_{i=1}p_i$ is the number of
vertices in the quiver of $A$. By Sturm's Theorem there are roots
$\lambda'_1<\lambda'_2<\cdots <\lambda'_{s-1}$ of
$q_{(p_1,\ldots,p_t-1)}(T)$ satisfying
$\lambda_1<\lambda'_1<\lambda_2<\cdots
<\lambda_{s-1}<\lambda'_{s-1}<\lambda_s$. The recursion formula
$$q_{(p_1,\ldots,p_t+1)}(T)=Tq_{(p_1,\ldots,p_t)}(T)-
q_{(p_1,\ldots,p_t-1)}(T)$$
implies that there is a root $\hat{\lambda}_i$ of
$q_{(p_1,\ldots,p_t+1)}(T)$ satisfying
$\lambda_i<\hat{\lambda}_i<\lambda_{i+1}$ for $i=1,\ldots,s-1$.

Indeed, the following picture illustrates the possible situations
concerning the distribution of roots\bigskip

\begin{figure}[h]
\begin{center}
\hskip6.6cm\includegraphics{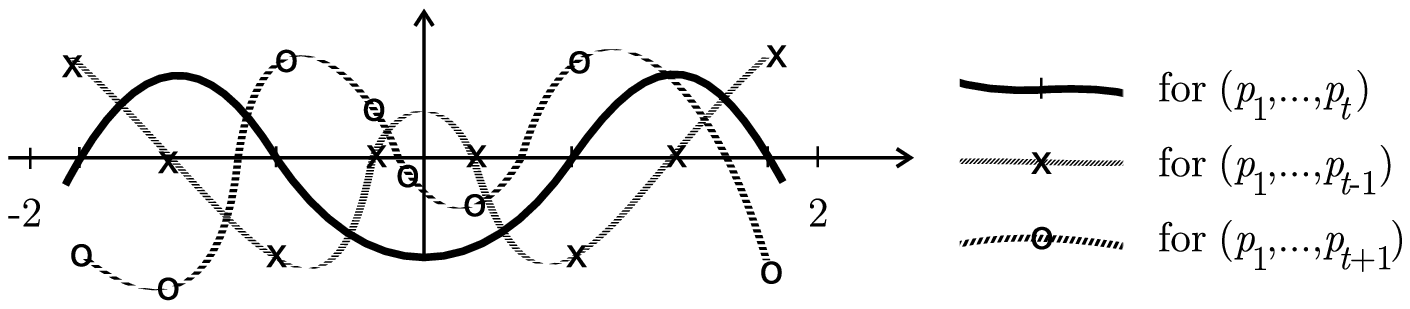}
\end{center}
\end{figure}
\bigskip
Accordingly, $\varphi_{\hat{A}}$ corresponding to the extended
canonical algebra $\hat{A}$ of type\break $(p_1,\ldots,p_t+1)$ has at
least $s-1\ge n-5$ eigenvalues in $\s^1$. By (1.4.b), at most $4$
roots of $\hat{f}_{(p_1,\ldots,p_t+1)}(T)$ are not in $\s^1$, as we
wanted to show.
\end{proof}

\end{document}